\documentclass[a4paper, 10pt]{amsart}

\usepackage{amsmath,amsthm,amssymb,graphicx,pinlabel,tikz,listings,xcolor}
\usepackage[all]{xy}
\usepackage{tikz-cd}

\SelectTips{cm}{}

\allowdisplaybreaks

\usepackage{hyperref}
\hypersetup{colorlinks=true,citecolor=specialblue,linkcolor=specialblue,urlcolor=black,filecolor=black}

\theoremstyle{definition}
\newtheorem{definition}{Definition}[section]
\newtheorem{remark}[definition]{Remark}

\theoremstyle{plain}

\newtheorem{proposition}[definition]{Proposition}
\newtheorem{lemma}[definition]{Lemma}

\makeatletter
    
    \@addtoreset{equation}{section}
  \makeatother

\makeatletter
\def\blfootnote{\xdef\@thefnmark{}\@footnotetext}

\definecolor{specialblue}{RGB}{35,88,180}
\definecolor{lightblue}{RGB}{10,180,255}
\definecolor{verylightblue}{RGB}{213,255,255}

\newcommand{\Z}{\mathbb{Z}}

\newcommand{\AJ}{\mathcal{A}J}
\newcommand{\calA}{\mathcal{A}}
\newcommand{\calI}{\mathcal{I}}

\newcommand{\calM}{\mathcal{M}}

\newcommand{\calG}{\mathcal{G}}

\newcommand{\calL}{\mathcal{L}}

\newcommand{\rAab}{r^{\calA, \mathfrak{a}}}

\newcommand{\Aut}{\operatorname{Aut}}

\newcommand{\Hom}{\operatorname{Hom}}
\newcommand{\Tr}{\operatorname{Tr}}
\newcommand{\Ker}{\operatorname{Ker}}

\newcommand{\Sp}{\operatorname{Sp}}

\newcommand{\GL}{\operatorname{GL}}

\author{Quentin Faes}

\begin{document}
\begin{abstract}
In this paper, we define trace-like operators on a subspace of the space of derivations of the free Lie algebra generated by the first homology group $H$ of a surface $\Sigma$. This definition depends on the choice of a Lagrangian of $H$, and we call these operators the \emph{Lagrangian traces}. We suppose that $\Sigma$ is the boundary of a handlebody with first homology group $H'$, and we show that the Lagrangian traces corresponding to the Lagrangian $\Ker (H \rightarrow H')$ vanish on the image by the Johnson homomorphisms of the elements of the Johnson filtration that extend to the handlebody.

\end{abstract}

\title{Lagrangian traces for the Johnson filtration of the handlebody group}
\maketitle

\vspace{-0.5cm}
 \tableofcontents
Let $\Sigma := \Sigma_{g,1}$ be a compact connected oriented surface of genus $g$ with one boundary component, and denote by $\calM:=\calM(\Sigma)$ the mapping class group of $\Sigma$. Let $V_g$ be a handlebody whose boundary minus a disk is identified with $\Sigma$. The handlebody group $\calA$ is the subgroup of mapping classes of $\calM$ that admit an extension to the handlebody. It is a non-normal subgroup of the mapping class group with infinite index, hence its study as a subgroup of $\calM$ is difficult. 

It is known, by the Dehn-Nielsen theorem, that $\calM$ can be studied via its action on the fundamental group $\pi := \pi_1(\Sigma_{g,1})$ of the suface, which is a free group on $2g$ generators. Let us denote $H := H_1(\Sigma_{g,1})$. Similarly, let $\pi' := \pi_1(V_g)$ and $H':= H_1(V_g)$ denote respectively the fundamental group and the first homology group of the handlebody. The abelians groups $H$ and $H'$ are the abelianization of $\pi$ and $\pi'$, respectively. We also denote $\mathbb{A}$ the kernel of the projection $\pi \rightarrow \pi'$ and $A$ the kernel of the projection $H \rightarrow H'$. We have the following commutative exact diagram, induced by the inclusion  $\Sigma_{g,1} \subset V_g$: 

\[ \begin{tikzcd} 0 \arrow[r]& \mathbb{A}\arrow[d] \arrow[r] &\pi \arrow[r] \arrow[d] & \pi' \arrow[r] \arrow[d , two heads] &0 \\0\arrow[r] & A\arrow[r] & H\arrow[r] & H'\arrow[r]& 0.  \end{tikzcd}\] The group $\calA$ actually coincides with the subgroup of $\calM$ that preserves $\mathbb{A}$ \cite{hensel}.

\begin{figure}[h]
	\centering
	\includegraphics[scale= 1]{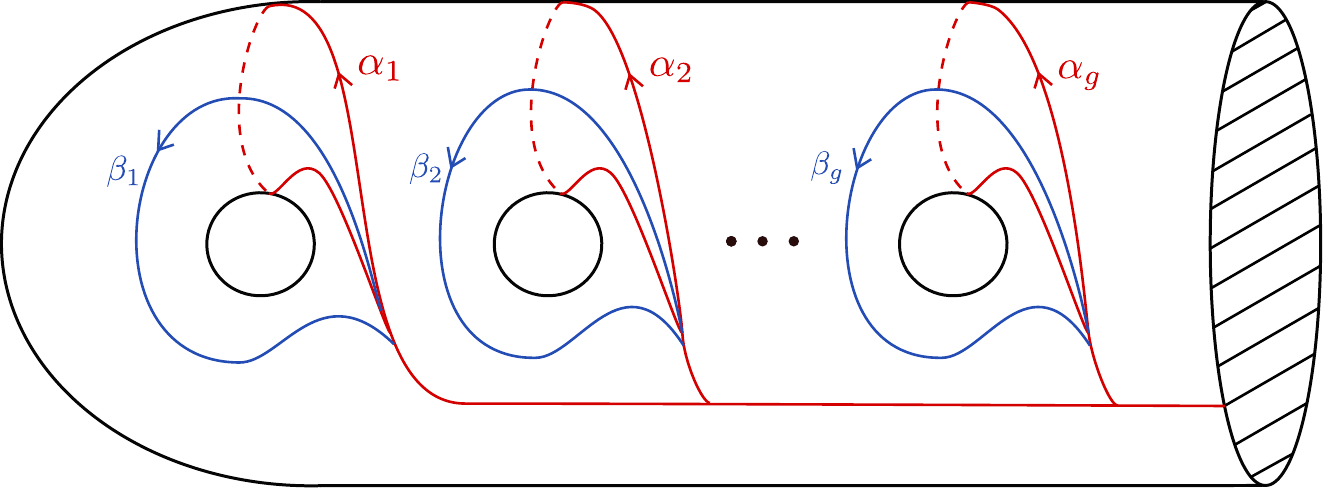}
	\caption{A basis of $ \pi$ adapted to the inclusion $\Sigma \subset V_g$.}
	\label{handlebody}
\end{figure}

Recall that the intersection form of $\Sigma$ induces a symplectic form $\omega : H \otimes H \rightarrow \Z$, hence $H$ is naturally isomorphic to $H^*$ via the map $x \mapsto \omega(x, - )$. The action of $\calM$ on the homology of the surface preserves the symplectic form $\omega$, inducing a representation $\rho_0 : \calM \rightarrow \Sp(H)$ whose kernel is, by definition, the Torelli group $\calI$. The duality $H \simeq H^*$ actually restricts to an isomorphism between the lagrangian $A$ and $H'^*$ (see \cite[Section 4.2]{faes}).
We fix a basis $(\alpha_1, \alpha_2 \dots \alpha_g, \beta_1, \beta_2 \dots \beta_g)$ of $\pi$, such that:
\begin{itemize}
\item it induces a symplectic basis of $H$.
\item $\mathbb{A}$ is normally generated by the $(\alpha_i)_{1 \leq i \leq g}$
\item the image $(\beta_i')_{1 \leq i \leq g}$ of the family $(\beta_i)_{1 \leq i \leq g}$ in $\pi'$ form a basis of $\pi'$.
\end{itemize}
\noindent Such a choice of basis for $\pi$ is induced by the curves drawn in Figure \ref{handlebody}. Notice that the curves $\alpha_i$ bound disks in the handlebody.

\blfootnote{This research has been supported by the project “AlMaRe” (ANR-19-CE40-0001-01) of the ANR and the project ``ITIQ-3D" of the Région Bourgogne Franche-Comté. The IMB receives support
from the EIPHI Graduate School (ANR-17-EURE-0002)}

The $k$-th term $J_k$ of the Johnson filtration is defined as the subgroup of $\calM$ of elements acting trivially on $\pi$ modulo the commutators of length $k+1$. In particular the first term of this filtration is $\calI$. The Johnson filtration $(J_k)_{k \geq 1}$ of the mapping class group is a separating filtration of $\calM$ and its associated graded space $\operatorname{Gr}^J := {\bigoplus}_{k\geq 1} \frac{J_k}{J_{k+1}}$ has been extensively studied. The latter is in particular a graded Lie ring. The so-called Johnson homomorphisms $(\tau_k)_{k \geq 1}$ induce an embedding $\tau$ of the graded space associated to the Johnson filtration to the graded Lie ring $D(H) := {\bigoplus}_{k\geq 1} D_k(H)$ of symplectic derivations of $\calL(H)$, the free graded Lie algebra generated by $H$. Recall that a derivation of the free Lie algebra can be identified via the restriction to $H$, with an element of $H^*  \otimes \calL (H) \simeq H \otimes \calL (H)$. The space $D_k(H)$ then corresponds to the kernel of the bracket from $H \otimes \calL_{k+1}(H)$  to $ \calL_{k+2}(H)$. The determination of the cokernel of $\tau$ is a question of major importance. In \cite{mor93}, Morita defined a family of trace-like operators that detect a part of this cokernel (but not all of it). 

In this paper, we are interested in the intersection of the Johnson filtration with $\calA$. We set $\calA_k := \calA \cap J_k$. The $\AJ$ \emph{filtration} $(\calA_k)_{k \geq 1}$ is also a separating filtration of $\mathcal{A}$. Hence we would like to understand its associated graded $\operatorname{Gr}^{\calA J}$, or, equivalently, the subspace $\tau_k(\calA_k) \subset D_k(H)$. Note that the study of the $\AJ$ filtration is important for the study of Heegaard splittings and finite-type invariants of integral homology 3-spheres \cite{pit,ric,faes2}. In \cite{omo}, Omori gave a generating set of the first term of the $\AJ$ filtration. In \cite{lev}, Levine noticed that the image of the restriction of $\tau_k$ to $\calA_k$ is included in the kernel $\calG_k$ of the natural projection $D_k(H) \rightarrow D_k(H')$, where $D_k(H') := \Ker \big( H' \otimes \calL_{k+1}(H') \rightarrow \calL_{k+2}(H')\big) $. Hence the Johnson homomorphisms induce a well-defined embedding $\operatorname{Gr}^{\calA J} \rightarrow {\bigoplus}_{k\geq 1} \calG_k$. A computation of Morita in \cite{mor} shows that $\calG_1 = \tau_1(\calA_1)$. We will show that this equality does not hold in general. In \cite{faes}, the author defined a family of trace-like operators $(\Tr^A_k)_{k \geq 1}$, defined on $\calG_k$, and valued in $S^k(H')$. We shall refer to these maps as the \emph{Lagrangian traces}. In the same paper, using a variation formula for the Casson invariant from \cite{mor}, the author shows that $\tau_2(\calA_2) = \Ker(\Tr^A_2) \cap \operatorname{Im}(\tau_2)$. Here, adapting the strategy of Morita in \cite{mor93} to the handlebody case, we prove the following theorem:\\

\noindent
\textbf{Theorem A.} \emph{Let $g \geq 2$, and $k \geq 2$, then the map $\Tr^A_k$ vanishes on $\tau_k(\calA_k) \subset \calG_k$.}\\

The proof of Theorem A uses the Magnus representation of the mapping class group, and also its reduction to the handlebody case. We will recall only what is needed here, and refer to \cite{saksurvey} for a survey on the subject. We define in Section \ref{sec2} two Magnus ``representations'' $r^A : \calA \rightarrow \operatorname{GL}(g, \mathbb{Z}[\pi'])$ and its abelian version $r^{A,\mathfrak{a}} : \calA \rightarrow \operatorname{GL}(g, \mathbb{Z}[H'])$. These are representations only in a twisted sense. For any group $G$, we shall denote $I_G$ the augmentation ideal of $\Z[G]$. For short, we set $I := I_\pi$. \\

\noindent
\textbf{Theorem B.} \emph{Let $g \geq 2$, and $\varphi \in \mathcal{A}_1$, then $\operatorname{det}(r^{A,\mathfrak{a}}(\varphi)) \in H' \subset \mathbb{Z}[H']$ and $$\operatorname{det}(r^{A,\mathfrak{a}}(\varphi)) = \Tr^A_1(\tau_1(\varphi)).$$ Consequently, the determinant of the Magnus representation is trivial on $\calA_k$, for $k \geq 2$.}\\

In the first section, we review the definition of the trace-like operators $(Tr^A_k)_{k \geq 1}$, and study their properties. In the second section, we will prove Theorem B and Theorem A by considering the magnus representation of the handlebody group, and relating it to the A-traces. Throughout the paper, we shall detail the connections with Morita's work in \cite{mor93}. \\

\noindent
\textbf{Acknowledgment.} The author would like to thank Gwénaël Massuyeau for his  corrections, as well as Renaud Detcherry and Jules Martel for discussions about the Magnus representation.

\section{The Lagrangian traces}
\label{sec1}
Before recalling the definition of the Lagrangian trace $\Tr^A_k$, let us recall a few facts about the Johnson homomorphisms $(\tau_k)_{k \geq 1}$ and the spaces $D_k(H)$. The $k$-th Johnson homomorphism of $\varphi \in J_k$ associates to the class $\{x\} \in H$ of an element $x$ in $\pi$, the class of $\phi(x)x^{-1}$ in $\Gamma_{k+1}\pi /\Gamma_{k+2} \pi \simeq \calL_{k+1}(H)$. This defines a map $ \tau_k : J_k \rightarrow \Hom(H,\calL_{k+1}(H))$. The target space can be identified with the space of derivations of degree $k$ of the free Lie algebra. Furthermore, one can prove that such a derivation sends $\omega$, the bivector dual to the symplectic form, to $0$. Such an element will be called a \emph{symplectic} derivation. These form a graded subalgebra $D(H)$. Any element in $D_k(H)$ has a \emph{Morita trace}, which coincides with the evaluation: 

\[ \begin{tikzcd}[cramped] D_k(H) \subset H^* \otimes \calL_{k+1}(H)\arrow[r] & H^* \otimes T_{k+1}(H)\arrow[r]& H^* \otimes H \otimes T_{k}(H) \arrow[r] & T_k(H) \arrow[r] & S^k(H) \end{tikzcd} \] where $\calL(H)$ embeds in the associative tensor algebra $T(H)$ via the rule $[a,b] = a\otimes b - b \otimes a$, and $H^* \otimes H \rightarrow \mathbb{Z}$ is the evaluation. Morita's trace map actually vanishes on the image of the Johnson homomorphism \cite{mor93}. Let us recall the definition of the trace-like operators defined in \cite{faes}. As we recalled, there is a canonical isomorphism between $A$ and $(H')^*$  given by the intersection pairing \[ \omega' : A \otimes H' \longrightarrow \mathbb{Z} \] \noindent which is induced by the symplectic form $\omega$. We denote $({\omega'})^{1,2}$ the contraction of the first two tensors in $A \otimes T(H')$ by $\omega'$. An element of $\mathcal{G}_k$ naturally induces an element of $A \otimes \mathcal{L}_{k+1}(H')$. Indeed $A \otimes \mathcal{L}_{k+1}(H') $ is the kernel of $H\otimes \mathcal{L}_{k+1}(H') \rightarrow H' \otimes \mathcal{L}_{k+1}(H')$. Then the map $\Tr^A_k$ is defined by 

\[
\begin{tikzcd}
 \operatorname{Tr}^A_k : \calG_k \arrow[r] &  A\otimes \mathcal{L}_{k+1}(H') \arrow[r, "i"]  & A \otimes T_{k+1}(H') \arrow[r, "(\omega')^{1,2}"]  &  T_{k}(H') \arrow[r, two heads] &  S^k(H').
 \end{tikzcd}\]

\begin{remark}
We actually have $A \simeq H'^*$, hence an element of $\calG_k$ naturally induces an element of $H'^* \otimes \calL_{k+1}(H')$, i.e. a derivation of $\calL(H')$. The Lagrangian traces that we define here are hence a handlebody version of the traces defined by Morita.
\label{remderiv}
\end{remark} 
Note that $\calG$ is a graded Lie subalgebra of $D(H)$. This subalgebra is not stable by the action of $\Sp(H)$, but it is stable by $\rho_0(\calA)$, the image of $\calA$ in $\Sp(H)$. This is precisely the subgroup of $\Sp(H)$ that preserves $A \subset H$, hence it acts on $H'$, and on $S^k(H')$.

\begin{lemma}
The linear maps $(\Tr_k^A)_{k \geq 1}$ are $\rho_0 (\calA)$-equivariant. 
\label{lemmaequiv}
\end{lemma}

\begin{proof}
It follows from the fact that an element of $\rho_0(\calA)$ preserves the intersection form $\omega'$.
\end{proof}

\noindent Note also that the graded space associated to the $\calA J$-filtration is stable by conjugation by elements of $\calA$, and the imbedding $\operatorname{Gr}^{\calA J} \rightarrow {\bigoplus}_{k\geq 1} \calG_k$ induced by $\tau$ is $\calA$-equivariant with respect to the actions we described.

For the sake of completeness, we describe how the Lagrangian traces behave with respect to the brackets in $\calG$. As a consequence of \cite[Lemma 3.2]{massak}, we have the following result:


\begin{proposition}
The map $\Tr^A$ vanishes on brackets of derivations. For any $\delta, \eta \in \calG$, we have \[ \Tr^A( [\delta, \eta]) = 0. \]
\end{proposition} \noindent We shall not use this proposition in the sequel.

\begin{remark}
Unlike the trace-like operators defined by Morita, the Lagrangian traces do not necessarily vanish in even degree. 
\end{remark}
\section{Magnus representations for the handlebody group}
\label{sec2}
\subsection{Morita's trace map and the Magnus representation}
We begin by reviewing the case of the mapping class group, following \cite{mor93}. In this section we work with use the notation $\gamma_i = \alpha_i$ and $\gamma_{g+i} =\beta_i$ for $ 1 \leq i \leq g$. Recall that the projections $(\beta_i')$ of the $\beta_is$ in $\pi'$ induce a basis of $\pi'$ and that the elements $\alpha_i$ normally generate the kernel of the projection $\pi \rightarrow \pi'$. We consider Fox derivatives $(\frac{\partial ~}{\partial \gamma_i}  : \Z[\pi] \rightarrow \Z[\pi])_{1 \leq i \leq g}$ with respect with the elements $\gamma_i$ of the basis. By definition, these are $\mathbb{Z}$-linear maps which verify $\frac{\partial \gamma_i}{\partial \gamma_j} = \delta_{i,j}$ for $1 \leq i,j \leq 2g$ and $\frac{\partial uv}{\partial \gamma_i} = \frac{\partial u}{\partial \gamma_i} + u\frac{\partial v}{\partial \gamma_i}$ for any $u,v \in \pi$. We also consider Fox derivatives for the basis chosen for $\pi'$.

To an element $\varphi$ in $\calM$, we associate the \emph{Fox matrix} $r(\varphi)$: 

\[\begin{array}{rl}
r :  \calM & \longrightarrow \GL(2g, \Z[\pi]) \\
\varphi &\longmapsto \bigg( \overline{\frac{\partial \varphi(\gamma_j)}{\partial \gamma_i}}\bigg)_{1 \leq i,j \leq 2g}
\end{array}\]
\noindent where for any $x \in \Z[\pi]$, $\overline{x}$ is the element obtained by applying the anti-automorphism induced by the inversion in $\pi$. The abelianization map from $\Z[\pi]$ to $\Z[H]$ then gives a map $$r^{\mathfrak{a}} :  \calM  \longrightarrow \GL(2g, \Z[H])$$usually called the \emph{Magnus representation} of the mapping class group. Nevertheless, this is \emph{not} a representation in the usual sense because it is a crossed homomorphism. When the map $r^{\mathfrak{a}}$ is restricted to the Torelli group $\mathcal{I}$, it becomes a representation. Recall that $D_1(H) \simeq \Lambda^3(H)$ and define $C$ to be the contraction map sending $a\wedge b \wedge c$ to $\omega(a,b)c + \omega(b,c) a + \omega(c,a)b$. Morita \cite[Prop. 6.15]{mor93} proved the following proposition: 

\begin{proposition}
Let $g \geq 2$, and $\varphi \in \calI$, then $\operatorname{det}(r^{\mathfrak{a}}(\varphi)) \in H \subset \mathbb{Z}[H]$ and $$\operatorname{det}(r^{\mathfrak{a}}(\varphi)) = 2C(\tau_1(\varphi)).$$ Consequently, the determinant of the Magnus representation is trivial on $J_k$, for $k \geq 2$.
\label{lemMor}
\end{proposition}

Furthermore, Morita proved that the reduction $r_k(\varphi)$ of the Fox matrix of $\varphi$ in $\GL(2g, \Z[\pi]/I^{k+1})$ can be recovered from the $k$-th Johnson homomorphism. Any element in $D_k(H)$ sends an element of $H$ to an element in $\Gamma_{k+1} \pi / \Gamma_{k+2} \pi$. The Fox derivatives of such an element relatively to the $\gamma_i's$ yield a well-defined vector in $\big(  I^{k} / I^{k+1}\big)^{2g}$. Indeed a Fox derivation sends an element of $\Gamma_k \pi$ to an element of $I^{k-1}$. For $d \in D_k(H)$, we set $$ \parallel d \parallel := \bigg(\frac{\partial d(\{\gamma_j\})}{\partial \gamma_i}  \bigg) \in \GL(2g, I^k/I^{k+1})$$ where $\{\gamma_j\}$ is the class of $\gamma_j$ in $H$. The following formula was proven by Morita:

\begin{equation}
\forall \varphi \in J_k, r_k(\varphi) = \operatorname{Id} + \overline{\parallel \tau_k(\varphi) \parallel}
\label{eqMor}
\end{equation}

\noindent For an element $d \in D_k(H)$, Morita's trace map was originally defined as the trace of $\parallel d \parallel$ after reduction in $I_H^k / I_H^{k+1}$ (the projection from $\pi$ to $H$ sends $I$ to $I_H$). From formula $\eqref{eqMor}$ it is clear that $$\operatorname{det}(r^{\mathfrak{a}}(\varphi)) \equiv 1 + \operatorname{Tr}(\overline{\parallel \tau_k(\varphi) \parallel}) \textit{ mod } I_H^{k+1}.$$ Then it follows from Proposition $\ref{lemMor}$ that Morita's trace map vanishes on the image of the Johnson homomorphisms. Note here that $\operatorname{Tr}$ is valued in $I_H^k / I_H^{k+1}$ which can be identified with $S^k(H)$.

\subsection{The handlebody case}
\label{sub:handle}

We will now adapt Morita's strategy of proof to the case of the handlebody group. We first give an equivalent definition for the A-trace $\Tr^A_k$. Let $d$ be an element of $\calG_k := \Ker(D_k(H) \rightarrow D_k(H'))$. As explained in Section \ref{sec1}, the class of $d$ in $H^* \otimes \mathcal{L}_{k+1}(H')$ sends $A$ to $0$. Equivalently, an element of $\calG_k$ induces an element of $H'^* \otimes \mathcal{L}_{k+1}(H')$. Similarly to the definition of $\parallel d \parallel$, we define $$ \parallel d \parallel^{A} := \bigg(\frac{\partial d(b_j')}{\partial \beta_i'}  \bigg) \in \GL(g, I_{\pi'}^k/I_{\pi'}^{k+1}).$$ It is clear that this is actually the reduction in $\pi'$ of the $g \times g$ right-down submatrix of $\parallel d \parallel$. We denote $\parallel d \parallel^{A,\mathfrak{a}} \in \GL(g, I_{H'}^k/I_{H'}^{k+1})$. Once again, note that $I_{H'}^k/I_{H'}^{k+1}$ can be identified with $S^k(H')$.

As a consequence \cite[Prop. 6.3]{mor93}, we have: 

\begin{lemma}
For any $k\geq 1$, for any $d \in \calG_k$, $\Tr^A_k(d) = \Tr(\parallel d \parallel^{A,\mathfrak{a}}) \in S^k(H').$
\label{lemmemorhan}
\end{lemma}

We now define the Magnus representation in the handlebody case. We are in fact interested in the right-down $g \times g$ part of the usual Fox matrix, after reduction in $\Z [\pi']$. Set \[\begin{array}{rl}
r^{\calA} :  \calA & \longrightarrow \GL(g, \Z[\pi']) \\
\varphi &\longmapsto \bigg( \overline{\frac{\partial \varphi(b_j')}{\partial b_i'}}\bigg)_{1 \leq i,j \leq 2g}
\end{array}\] and $r^{\calA, \mathfrak{a}} : \calA \rightarrow \GL(g, \Z[H'])$ its abelian version. It is actually the composition of $\calA \rightarrow \Aut (\pi')$ with the usual Magnus representation of the group $\Aut (\pi')$. The map $r^\calA$ is a crossed homomorphism in the following sense:

\begin{equation}  
\forall \varphi, \psi \in \calA, r^{\calA}(\varphi \psi) = r^{\calA}(\varphi) \big(\varphi \cdot r^{\calA}(\psi)\big)
\label{crossed}
\end{equation} where $\varphi \cdot r^{\calA}(\psi)$ is the matrix obtained from $r^{\calA}(\psi)$ after applying $\varphi$ to all the coefficients.

\begin{remark}
The mapping $\GL(2g, \Z[\pi]) \rightarrow \GL(g, \Z[\pi'])$ (resp. its abelian version), which consist in projecting in $\pi'$ (resp. $H'$) the coefficients of the righ-down matrix is of course not a homomorphism. This explains why the Lagrangian traces are not reductions of the Morita trace maps (see \cite[Remark 4.4]{faes}).
\end{remark} 

We are now ready to prove Theorem B.

\begin{figure}[h]
	\centering
	\includegraphics[scale= 1]{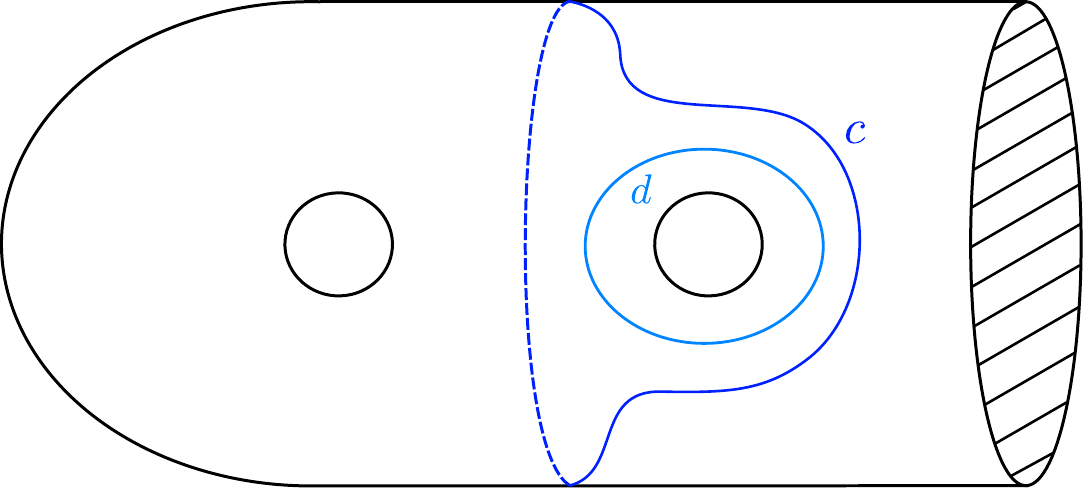}
	\caption{Two curves $c$ and $d$ such that $T_c\circ T_d^{-1}$ normally generates $\calA_1$.}
	\label{lagtr}
\end{figure}

\begin{proof}[Proof of Theorem B]
First, the map $ \operatorname{det} \circ \rAab$ is $\calA$-equivariant. Indeed, if $\psi$ is an element of $\calA$, and $\varphi$ is at least in $\calA_1$, we deduce from \eqref{crossed} that \[r^{\calA, \mathfrak{a}}(\psi\varphi \psi^{-1}) = r^{\calA, \mathfrak{a}}(\psi) \big(\psi \cdot r^{\calA, \mathfrak{a}}(\varphi)\big) \big(\psi \cdot r^{\calA, \mathfrak{a}}(\psi^{-1})\big),\] because $\varphi$ acts trivially on $H'$. We conclude by observing that \[\operatorname{Id} = r^{\calA, \mathfrak{a}}(\psi \psi^{-1}) = r^{\calA, \mathfrak{a}}(\psi) \big(\psi \cdot r^{\calA, \mathfrak{a}}(\psi^{-1})\big),\] hence $r^{\calA, \mathfrak{a}}(\psi\varphi \psi^{-1})$ is conjugated to $\psi \cdot r^{\calA, \mathfrak{a}}(\varphi)$. The map $\Tr^A_1 \circ \tau_1$ is also $\calA$-equivariant, by Lemma \ref{lemmaequiv}. According to \cite{omo}, the subgroup $\calA_1$ is normally generated by the map $\varphi := T_c \circ T_d^{-1}$ described in Figure \ref{lagtr}. Note that $\varphi$ is an \emph{annulus twist} in the sense of \cite{hensel} and a BP map in the sense  of \cite{joh79}. Indeed $c$ and $d$ cobound an annulus in the handlebody, and a surface in $\Sigma$, which ensures that $\varphi$ is both an element of $\calA$ and $\calI$. We set $\tilde{c} := \beta_2[\alpha_1, \beta_1]$ a lift of the curve $c$ in $\pi$. The action of $\varphi$ on $\pi$ is
\[\begin{array}{rl}
\alpha_1 & \mapsto \tilde{c}\alpha_1 \tilde{c}^{-1}\\
\alpha_2 & \mapsto \alpha_2 \beta_2\tilde{c}^{-1} \\
\beta_1 & \mapsto \tilde{c} \beta_1 \tilde{c}^{-1}\\
\beta_2 & \mapsto \tilde{c} \beta_2 \tilde{c}^{-1},\\
\end{array}\] hence we compute directly that $$ \rAab(\varphi) = \begin{pmatrix}
b_2'^{-1} & 0\\
1-b_1'^{-1} & 1
\end{pmatrix}.$$ We deduce that $\operatorname{det}(\rAab(\varphi)) = b_2'^{-1} \in \Z[H']$, and, switching to the additive notation in $H'$, we obtain an element $-b_2' \in H'$. By equivariance, as the map $\operatorname{det} \circ \rAab$ is a multiplicative map, we deduce that it sends the subgroup $\calA_1$ in $H'$. Furthermore, by \cite{joh80}, we know that $\tau_1(\varphi) = a_1 \wedge b_1 \wedge b_2$. Hence, $\Tr^A_1( \tau_1(\varphi)) = - b_2'$, and we deduce that, on $\calA_1$, $$\operatorname{det} \circ r^{A,\mathfrak{a}} = \Tr^A_1 \circ \tau_1.$$ As the first Johnson homomorphism vanishes on $\calA_k$, for $k \geq 2$, we deduce that the map $\operatorname{det} \circ r^{\calA,\mathfrak{a}}$ is trivial on these subgroups.
\end{proof}

We can now prove Theorem A.

\begin{proof}[Proof of Theorem A]
The argument is similar to \cite{mor93}. Indeed, we set $r_k^\calA$ the reduction of $r^\calA$ modulo $I_{\pi'}^{k+1}$. Clearly, as in \cite{mor93}[Prop. 6.13], we have the following equality: \begin{equation}
\forall \varphi \in \calA_k, r_k^{\calA}(\varphi) = \operatorname{Id} + \overline{\parallel \tau_k(\varphi) \parallel^A}.
\label{eqMorhan}
\end{equation}
\noindent Recall that $\varphi \in \calA_k$ implies that $\tau_k(\varphi) \in \calG_k$, hence the right-hand side of the equation makes sense.  We can deduce equation \eqref{eqMorhan} from equation \eqref{eqMor}, by reducing the coefficients of the right-down $g\times g$ matrix of $r_k(\varphi)$, or directly using the very definition of the Johnson homomorphisms. After projection of the coefficients in $\Z[H']$, we obtain $$\operatorname{det}(r^{\mathfrak{a}, \calA}(\varphi)) \equiv 1 + \operatorname{Tr}(\overline{\parallel \tau_k(\varphi) \parallel^A}) \textit{ mod } I_{H'}^{k+1}.$$ Using Theorem B and Lemma \ref{lemmemorhan}, we infer that, for any $\varphi \in \calA_k,$ $\Tr^A_k(\tau_k(\varphi)) \equiv 0 \textit{ mod } I^{k+1}_{H'}$.  

\end{proof}

\bibliographystyle{plain}
\bibliography{manuscrit}
\end{document}